\documentclass{amsart}
\usepackage[all]{xy}
\usepackage{amsmath}
\usepackage{amsfonts,amssymb}
\usepackage{latexsym}
\usepackage{pb-diagram}
\usepackage{url}
\RequirePackage[T1]{fontenc}

\newcommand{\al}{{\alpha}}

\newcommand{\io}{{\iota}}

\newcommand{\brd}{\bar{\pa}}

\newcommand{\DGL}{{\bf DG}\textbf{-}{\bf Lie}}
\newcommand{\Set}{{\bf Sets}}

\newcommand{\Nilp}{{\bf Nilp}}

\newtheorem{dfn}{Definition}[section]
\newtheorem{prop}[dfn]{Proposition}
\newtheorem{cor}[dfn]{Corollary}

\newtheorem{thm}[dfn]{Theorem}
\newtheorem{lem}[dfn]{Lemma}
\newtheorem{pb}[dfn]{Problem}

\newtheorem{rmq}[dfn]{Remark}
\newtheorem{ex}[dfn]{Example}

\newcommand{\bul}{{\bullet}}

\newcommand{\ad}{\mathop{\rm{ad}}}

\newcommand{\R}{\mathbb{R}}

\newcommand{\N}{\mathbb{N}}

\newcommand{\C}{\mathbb{C}}

\newcommand{\cL}{\mathcal{L}}
\newcommand{\cM}{\mathcal{M}}
\newcommand{\cC}{\mathcal{C}}

\newcommand{\cA}{\mathcal{A}}
\newcommand{\mO}{\mathcal O}

\newcommand{\Edif}{~\!\!^\cL\!d}

\newcommand{\Elie}{~\!\!^\cL\!L}
\newcommand{\EOm}{~\!\!^\cL\!\Omega}

\newcommand{\EA}{~\!\!^\cL\!\cA}

\newcommand{\pa}{\partial}
\newcommand{\mg}{\mathfrak g}
\newcommand{\mh}{\mathfrak h}
\newcommand{\mm}{\mathfrak m}

\author{Damien Calaque and Gilles Halbout}
\thanks{D.C.~is on leave of absence from Institut Camille Jordan UMR5208, Universit\'e Lyon1, F-69622 Villeurbanne, France}

\title{Weak quantization of Poisson structures}

\begin{document}

\begin{abstract}
In this paper we prove that any Poisson structure on a sheaf of Lie algebroids admits a weak 
deformation quantization, and give a sufficient condition for such a Poisson structure to 
admit an actual deformation quantization. 
We also answer the corresponding classification problems. In the complex symplectic case, 
we recover in particular some results of Nest-Tsygan and Polesello-Schapira. 

We begin the paper with a recollection of known facts about deformation theory of 
cosimplicial differential graded Lie algebras. 
\end{abstract}

\maketitle

\tableofcontents

\section*{Introduction}

In this paper we prove a very general result concerning the deformation 
quantization problem for sheaves of Lie algebroids.

Following \cite{Hin}, any reasonable formal deformation problem can be described by a 
functor on differential graded (DG) artinian rings with values in simplicial sets, 
representable by some DG Lie (or perhaps $L_\infty$) algebra. 
In this paper we deal with a deformation problem that is described by a {\it sheaf of} 
differential graded Lie algebras (DGLA). We solve this problem and take this opportunity to
recall the construction of the deformation functor associated to cosimplicial DGLA. 

Deformation quantization problem for a $C^\infty$ Poisson manifold $(M,\pi)$ has been solved 
by Kontsevich in \cite{K1}. Kontsevich first proves a formula in the local 
case ($M=\R^d$) and then apply an appropriate globalization procedure. Actually the existence 
of a deformation quantization is a part of a more general picture: Kontsevich proves in 
\cite{K1} that the DGLA of poly-differential operators is formal. 

This formality theorem is generalized to a large class of sheaves of Lie algebroids in 
\cite{CVdB} (see also \cite{C,C2,CDH} for the particular cases of $C^\infty$ and holomorphic 
Lie algebroids). In this paper, we prove that any Poisson structure admits a weak deformation 
quantization (Theorem \ref{thm-main}). We also give a sufficient condition for such a Poisson 
structure to admit an actual deformation quantization. 
We also answer the corresponding classification problems. In the complex symplectic case, 
we recover in particular some results of Nest-Tsygan and Polesello-Schapira. \\

This paper can be seen as an attempt to understand some claims of \cite{K2} where this question is discussed 
in the context of algebraic geometry. 
We also want to emphazise the great importance of the extremely enlighting ``homotopical point-of-view'' 
\cite{Hin} on deformation theory. 

Throughtout the paper $k$ is a field with $char(k)=0$. \\

Plan of the paper. 
In section 1 we review some basic materials concerning models for (cosimplicial) simplicial sets and 
(cosimplicial) DG Lie algebras. We also define the deformation functor associated to a cosimplicial DG 
Lie algebra. 

In section 2 we recall the construction of the Deligne 2-groupoid associated to ``quantum type'' DG Lie algebras 
and give an explicit description of the deformation functor associated in this situation.

In section 3 we apply the previous constructions to some quantization problems. Namely, we first prove that any 
Poisson structure associated to a (locally free of finite rank) Lie algebroid admits a weak quantization. 
We then classify such weak quantizations and then give a sufficient condition for the existence of an usual 
quantization of a given Poisson structure. We compare our results with previous works \cite{K2,Po,PS,NT}. 

\subsection*{Acknowledgements}

Both authors are grateful to their former host institution, IRMA (Strasbourg), where they 
started this project. \\
D.C.~heartly thanks Mathieu Anel for teaching him model categories and modern homotopy theory. He is indebted 
to Amnon Yekutieli for reference \cite{LY} and many fruitful e-mail discussions. He also thanks ETH (Z\"urich) 
and IHES (Bures-sur-Yvette), where part of this work was improved, for hospitality. His work has been partially 
supported by the European Union through the FP6 Marie Curie RTN ENIGMA (Contract number MRTN-CT-2004-5652). \\
The authors also thank Vasiliy Dolgushev for reference \cite{BGNT}. 

\section{Basic materials}

For the reader who wants to learn about model categories we refer to the very down-to-earth introduction 
\cite{GK} and references therein. 

\subsection{Model categories and (co)simplicial methods}

Let us first recall that a {\it closed model category (CMC)} is a category equipped with three classes of 
morphisms (called {\it fibrations}, {\it cofibrations}, and {\it weak equivalences}) satisfying the axioms 
(CM1)-(CM5) of \cite{Q2}. 

\subsubsection{(Co)Simplicial objects}

Let $\cC$ be a category. Let us denote by $c\cC=\cC^{{\bf\Delta}}$ (resp.~$s\cC=\cC^{{\bf\Delta}^{{\rm op}}}$) 
the category of cosimplicial (resp.~simplicial) objects in $\cC$. 
Here ${\bf\Delta}$ denotes the {\it ordinal number category} (or {\it simplicial category}), i.e.~the 
category with objects ordered finite sets $[k]=\{0,\dots,k\}$ and morphisms (weakly) order 
preserving maps. In other words, 
$$
{\bf\Delta}_k^l:={\rm Hom}_{{\bf\Delta}}([k],[l])=\{(i_1,\dots,i_k)|0\leq i_1\leq\dots\leq i_k\leq l\}\,.
$$

Among all morphisms in ${\bf\Delta}$ there are the following remarkable ones: for $i\in[k]$, 
$\delta_i=(0,\dots,i-1,i+1,\dots,k):[k-1]\to[k]$ and $\sigma_i=(0,\dots,i,i\dots,k):[k+1]\to[k]$. 
Moreover, any morphism in ${\bf\Delta}$ 
is a composition of these two types of morphisms. The following easy lemma provides a full list of 
relations between them. 
\begin{lem}
One has the following identities in ${\bf\Delta}$: 
\begin{itemize}
\item $\delta_i\delta_j=\delta_j\delta_{i-1}$ if $i>j$, 
\item $\sigma_i\sigma_j=\sigma_j\sigma_{i+1}$ if $i\geq j$, 
\item  $\sigma_i\delta_j=
\begin{cases}
\delta_{j-1}\sigma_i\quad(\textrm{ if }i<j-1) \\
{\rm id}\quad(\textrm{ if }i=j-1,j) \\
\delta_j\sigma_{i-1}\quad(\textrm{ if }i>j)
\end{cases}$
\end{itemize}
\end{lem}

As a matter of notation, for cosimplicial (resp. simplicial) objects $C$, we will write 
$C^k:=C([k])$ (resp.~$C_k:=C([k])$) for $k\geq0$, and $f:=C(f)$ (resp.~$f^*:=C(f)$) for any morphism 
$f$ in ${\bf\Delta}$. 

\begin{ex}\emph{
For any $n\in\N$ we define the {\it geometric $n$-simplex} $\sigma^n$ as the convex hull of the canonical basis 
in $\R^{n+1}$, and we identify $[n]$ with the vertices of $\sigma^n$. Then any map $\phi:[m]\to[n]$ can be 
extended to a PL map $\phi:\sigma^m\to\sigma^n$. Thus $(\sigma^n)_n$ is a cosimplicial PL manifold. }
\end{ex}

\subsubsection{The model category of simplicial sets}

Let us denote by $\Set$ the category of sets. The category $s\Set$ of simplicial sets 
has a natural structure of a closed model category that we are going to describe. 

For this we need to introduce some remarkable objects in $s\Set$. 
The {\it standard $n$-simplex} is the simplicial set 
$\Delta^n={\rm Hom}_{{\bf\Delta}}(-,[n]):{\bf\Delta}^{{\rm op}}\to\Set$. 
In other words, $\Delta^n_k={{\bf\Delta}}_k^n$. Moreover any $f\in{\bf\Delta}_m^n$ induces a 
morphism of simplicial sets $\Delta^m\to\Delta^n$. Therefore one can also consider the 
{\it boundary} $\pa\Delta^n:=\cup_{k=0}^n\delta_k\Delta^{n-1}\subset\Delta^n$ of 
$\Delta^n$, and its {\it horns} 
$\Lambda^{i,n}:=\cup_{k\neq i}\delta_k\Delta^{n-1}\subset\Delta^n$ ($i\in[n]$). 
In particular $\Delta^n_0=(\partial\Delta^n)_0={[n]}$ for $n\neq0$. 

For any simplicial set $X_\bullet$ one can define define its {\it realization} $|X|$ as a certain 
colimit in a category of topological spaces (see \cite[Definition 1.19]{GK} ; simply note that 
$|\Delta^n|=\sigma^n$). Then one can define the set of path components $\pi_0(X):=\pi_0(|X|)$ and 
homotopy groups $\pi_i(X,x):=\pi_i(|X|,x)$ ($i>0$) for any $x\in\pi_0(X)$. Any morphism of simplicial 
sets $f:X\to Y$ induces morphisms $f_0:\pi_0(X)\to\pi_0(Y)$ and $\pi_i(X,x)\to\pi_i(Y,f_0(x))$ for $i>0$. 

The CMC structure on $s\Set$ is such that the morphism $f:X\to Y$ is a 
\begin{itemize}
\item weak equivalence if it induces isomorphisms $\pi_0(X)\cong\pi_0(Y)$ and 
$\pi_n(X,x)\cong\pi_n(Y,f_0(x))$ for any $x\in\pi_0(X)$ and any $n>0$, 
\item cofibration if $f_n$ is injective for any $n\geq1$. 
\end{itemize}

A fibrant object $X$ in $s\Set$ is called a {\it weak $\infty$-groupoid} (or {\em Kan complex}). It has 
the property that $\pi_0(X)=X_0/\sim$, where $x\sim y$ if and only if there exists $\gamma\in X_1$ such 
that $\delta_1^*\gamma=x$ and $\delta_0^*\gamma=y$. 

Moreover, a simplicial set $C$ is fibrant if and only if the following condition is met: 
\begin{quote}
{\it every morphism $\Lambda^{i,n}\to C$ can be extended to a morphism $\Delta^n\to C$. }
\end{quote}

\medskip

We will also need the notion of {\it simplicial closed model category (SCMC)}: it is a CMC category $\cM$ 
enriched over simplicial sets (we denote by ${\rm Hom}_\cM^\Delta(X,Y)$ the enriched Hom space) that satifies 
Quillen's axiom (SM7) of \cite{Q1}. The category of simplicial sets is naturally a SCMC, where the CMC structure 
is the previous one and the simplicial structure is given by 
$$
{\rm Hom}_{s\Set}^\Delta(C,D)_n:={\rm Hom}_{s\Set}(\Delta^n\times C,D)\,.
$$

\subsubsection{Reedy model categories}\label{sec-Reedy}

The result of this paragraph first appeared in the unpublished paper \cite{Ree}. 
Let $\cM$ be a CMC. We want to describe a model structure on $c\cM$ 
(in his paper Reedy deals with $s\cM=(c(\cM^{\rm op}))^{\rm op}$). 

First of all, for any object $C$ in $c\cM$ we define its {\em $n$-th matching object} in $\cM$ to be the colimit 
$M^nC:=\underrightarrow{\lim}C^k$ taken over all surjections $[n]\to[k]$ in ${\bf\Delta}_k^n$ with $k<n$. 
We have natural morphisms $C^n\to M^nC$. 

Then $c\cM$ has a CMC structure with a morphism $f:C\to D$ being a 
\begin{itemize}
\item weak equivalence if $f^n:C^n\to D^n$ is a weak equivalence in $\cM$ for all $n\geq0$, 
\item fibration if the induced map $C^n\to D^n\times_{M^nD}M^nC$ is a fibration in $\cM$ for all $n\geq0$. 
\end{itemize}
\begin{rmq}
One can of course define the {\em $n$-th latching object} of $C$ to be the limit $L^nC:=\underleftarrow{\lim}C^k$ 
taken over all injections $[k]\to[n]$ in ${\bf\Delta}_n^k$ with $k<n$, and then cofibrations are characterized in a 
similar way as fibrations. 
\end{rmq}

A fibrant object $C$ in $c\cM$ is an object such that $C^n\to M^nC$ is a fibration for all $n$ (this is 
an abstract way of describing descent condition). Moreover it is a standard fact that a fibrant object 
$C$ in $c\cM$ is termwize fibrant, that is to say $C^n$ is fibrant for any $n\geq0$ (the converse 
is obviously false). 

\subsubsection{Cosimplicial simplicial sets}

It follows from the previous two paragraphs that the category $cs\Set$ of cosimplicial simplicial sets has 
a natural structure of a CMC. Notice that there is a remarkable cofibrant object $\Delta$ in 
$cs\Set$ defined as follows: 
\begin{eqnarray*}
\Delta\,:\,{\bf\Delta} & \longrightarrow & s\Set \\
{[n]} & \longmapsto & \Delta^n={\rm Hom}_{{\bf\Delta}}(-,[n]) \\
f\in{\bf\Delta}_m^n & \longmapsto & (g\in\Delta_k^m\mapsto f\circ g\in\Delta_k^n)\,.
\end{eqnarray*}
Using $\Delta$ we define the total space functor ${\rm Tot}:cs\Set\longrightarrow s\Set$ to be 
$$
C\longmapsto {\rm Tot}(C):={\rm Hom}^\Delta_{cs\Set}(\Delta,C)\,. 
$$
The functor ${\rm Tot}$ preserves fibrations, trivial fibrations and weak equivalences 
between fibrant objects (as $\Delta$ is cofibrant). 

The set ${\rm Tot}(C)_0$ of $0$-simplices is given by sequences $(\alpha_0,\alpha_1,\dots,\alpha_n,\dots)$ 
such that $\alpha_n\in C_n^n$ and $s\alpha_i=s^*\alpha_j$ (in $C_i^j$) for any $s\in{\bf\Delta}_i^j$ 
(this is ${\rm Hom}_{cs\Set}(\Delta,C)$). 

\subsection{Homotopy theory of DG Lie algebras}

\subsubsection{The model category of DG Lie algebras}

Let us denote by $\DGL$ the category of DG Lie algebras with morphisms being standard morphisms 
of DG Lie algebras. $\DGL$ admits a closed model structure with a morphism $f:\mg\to\mh$  being a 
\begin{itemize}
\item weak equivalence if $H(f)$ is an isomorphism\footnote{In other words, weak equivalences are 
quasi-isomorphisms. }, 
\item fibration if $f$ is surjective. 
\end{itemize}
In particular one can see that any object in $\DGL$ is fibrant. 

Moreover $\DGL$ can be enriched over simplicial sets so that it becomes a SCMC. Namely, 
$$
{\rm Hom}^\Delta_{\DGL}(\mg,\mh)_n:={\rm Hom}_{\DGL}(\mg,\Omega_n\otimes\mh)\,,
$$
with $\Omega_n$ being the DG commutative algebra of differential forms on the geometric $n$-simplex (in particular 
$(\Omega_n)_n$ is a simplicial DG commutative algebra). 

%
%
%

\subsubsection{Hinich's deformation functor}

Let us denote by $\Nilp$ the category of pronilpotent commutative $k$-algebras, and by $s\Set^\Nilp$ the SCMC of 
functors $\Nilp\to s\Set$: weak equivalences (resp.~(co)fibrations) in $s\Set^\Nilp$ are termwize weak equivalences 
(resp.~(co)fibrations). 

Following \cite{Hin} we define a functor 
$$
\DGL\longrightarrow s\Set^\Nilp\,;\,\mg\longmapsto\Sigma_\mg
$$
that preserves fibrations and weak equivalence. Namely, for any DG Lie algebra $(\mg,d,\mu)$ and any pronilpotent 
commutative algebra $\mm$ the set $\Sigma_\mg(\mm)_n$ is the set of degree one elements $\Pi$ in 
$\Omega_n\otimes\mg\otimes\mm$ such that 
$$
d\Pi+\frac12\mu(\Pi,\Pi)=0\,.
$$
In other words, $\Sigma_\mg(\mm)_n$ is the set of {\it Maurer-Cartan elements} of the DG Lie algebra 
$\Omega_n\otimes\mg\otimes\mm$. 

\subsubsection{The Maurer-Cartan functor associated to a cosimplicial DGLA}

We consider the category $c\DGL$ of cosimplicial DG Lie algebras with its Reedy model structure. $\Sigma$ 
obviously extends to a functor $c\DGL\to cs\Set^\Nilp$ preserving fibrations and weak equivalences. Composing it 
with ${\rm Tot}$ we obtain a functor 
$$
Def\,:\,c\DGL\longrightarrow s\Set^\Nilp\,;\,Def_\mg(\mm):={\rm Tot}(\Sigma_\mg(\mm))
$$
that preserves fibrations, trivial fibrations and weak equivalences between fibrant objects. 
We finally define a functor 
$$
\underline{MC}:=\pi_0\circ Def\,:\,c\DGL\longrightarrow\Set^\Nilp\,,
$$
that sends weak equivalences between fibrant objects to equalities. We call $\underline{MC}_\mg$ the 
{\it Maurer-Cartan functor} of $\mg$, it is an homotopy invariant. 
\begin{rmq}\label{rmq-sheaf}
Let $X$ be a topological space and $X=\bigcup_{i\in I}U_i$ an open cover  indexed by a totally order set $I$. 
Then for any presheaf of DG Lie algebras $U\mapsto\mg(U)$ one can construct a cosimplicial DG Lie algebra 
$\mg(\underline{U}^\bul)$ in an obvious way. Moreover if $U\mapsto\mg(U)$ is a sheaf then 
$\mg(\underline{U}^\bul)$ is fibrant (conversely, if $\mg(\underline{U}^\bul)$ is fibrant for ANY open cover 
then $U\mapsto\mg(U)$ is a sheaf). 
\end{rmq}

\subsubsection{Sheaves of DG Lie algebras}

Given a site $\mathcal C$ one can define a model structure on the category $\mathcal P_{\mathcal C}(\DGL)$ of 
{\it presheaves of DG Lie algebras} over $\mathcal C$ (see \cite{Hin2}) in such a way that fibrant objects are 
precisely sheaves of DG Lie algebras. One could repeat the previous constructions in this context. 

Nevertheless, if the site $\mathcal C$ is not too bad (e.g.~if it is the small site of a reasonnable topological 
space) then this description is more or less equivalent to the cosimplicial one (thanks to remark \ref{rmq-sheaf}). 

\section{Quantum type DG Lie algebras and $2$-groupoids}

Following the terminology of \cite{LY}, by a {\it quantum type DG Lie algebra} $\mg$ we mean a 
DG Lie algebra such that $\mg^{[i]}=0$ for $i<-1$. Getzler shows in \cite{Ge2} that if $\mg$ is a 
quantum type DG Lie algebra then $\Sigma_\mg(\mm)$ is weakly equivalent to the nerve of a strict $2$-groupoid, 
called the {\em Deligne $2$-groupoid} of $\mg\otimes\mm$ (see \cite{Ge1}). In this section we recall the 
construction of this $2$-groupoid and then give an explicit description of the functor $\underline{MC}$ 
in the quantum type situation. 

This description already appears in a slightly different formulation in \cite[Section 3]{BGNT} (in particular one 
can recover Theorem 3.6 of \cite{BGNT} from our approach). 

\subsection{The Deligne $2$-groupoid of a quantum type DGLA}

We follow \cite{Ge1}. 
Let $(\mg,d,[,])$ be a DG Lie algebra and $\mm$ a pronilpotent commutative $k$-algebra. We are going to define 
a $2$-groupoid $Del_\mg(\mm)$. 
Objects in $Del_\mg(\mm)$ are given by the set of Maurer-Cartan elements, i.e.~ elements 
$\Pi\in(\mg\otimes\mm)^{[1]}$ such that 
$$
d\Pi+\frac12[\Pi,\Pi]=0\,.
$$
The prounipotent group $\exp((\mg\otimes\mm)^{[0]})$ acts on objects in the following way: 
for any $q\in(\mg\otimes\mm)^{[0]}$ and any $\Pi\in(\mg\otimes\mm)^{[1]}$ 
$$
\exp(q)\cdot\Pi:=\Pi -\sum_{n=0}^\infty\frac{{\ad(q)}^n}{(n+1)!}d_\Pi q\,,
$$
where $d_\Pi(q)=d q +[\Pi,q]$.
The subset of Maurer-Cartan elements is obviously stable under this action\footnote{This action is the 
exponentiation of the infinitesimal affine action $q\cdot\Pi=dq+[q,\Pi]$. }. The translation groupoid 
associated to the group action $\exp((\mg\otimes\mm)^{[0]}) \times MC(\mg,\mm) \to MC(\mg,\mm)$ is the underlying 
$1$-groupoid of Deligne $2$-groupoid in the following sense: the set of morphisms between $\Pi$ and 
$\exp(q)\cdot\Pi$ of the structure of a $1$-groupoid, which is the translation groupoid associated to 
the group action 
$$
\exp((\mg\otimes\mm)_\Pi^{[-1]})\times\exp((\mg\otimes\mm)^{[0]}) \to \exp((\mg\otimes\mm)^{[0]})\,.
$$
Here $(\mg\otimes\mm)_\Pi^{[-1]}$ denotes the Lie algebra $((\mg\otimes\mm)^{[-1]},[-,-]_\Pi)$, where 
$[u,v]_\Pi=[d_\Pi u,v]$, and the group action 
$$
\exp(u)\cdot\exp(q):=\exp(q)\exp(d_\Pi u)\,.
$$
comes from the Lie algebra morphism $d_\Pi:(\mg\otimes\mm)_\Pi^{[-1]}\to(\mg\otimes\mm)^{[0]}$. 

For $\Pi\in MC(\mg,\mm)$, $q\in(\mg\otimes\mm)^{[0]}$ and $u,v\in(\mg\otimes\mm)_\Pi^{[-1]}$, the 
{\it vertical} composition of $2$-morphisms is given by the formula 
$$
(\exp(v),\exp(q)\exp(d_\Pi u),\Pi)\circ_{\textrm{v}}(\exp(u),\exp(q),\Pi)
=(\exp(u)\exp(v),\exp(q),\Pi)\,;
$$
$$
\xymatrix@C=13pt@R=-2pt@M=6pt{
& \ar@{=>}[ddd]^{e^u} &&&&& \ar@{=>}[dddddd]^{e^ue^v} & \\
&&&&&&& \\
&&&&&&& \\
\Pi\ar@/^3pc/ [rrr]^{e^q} \ar[rrr]^{e^qe^{d_\Pi u}} \ar@/_3pc/ [rrr]_{e^qe^{d_\Pi u}e^{d_\Pi v}} & 
\ar@{=>}[ddd]^{e^v} && e^q\cdot\Pi \ar@{~>}[rr] & & 
\Pi\ar@/^3pc/ [rr]^{e^q} \ar@/_3pc/ [rr]_{e^qe^{d_\Pi u}e^{d_\Pi v}} && e^q\cdot\Pi \\ 
&&&&&&& \\
&&&&&&& \\
&&&&&&&}
$$
$$
\texttt{Figure 1.~Vertical composition}
$$
~\\
and for $q'\in(\mg\otimes\mm)^{[0]}$ and $u'\in(\mg\otimes\mm)_{\exp(q)\cdot\Pi}^{[-1]}$ the {\it horizontal} 
composition of $2$-morphisms is given by 
\begin{align*}
(\exp(u'),\exp(q'),\exp(q).\Pi)\circ_{\textrm{h}}(\exp(u),\exp(q),\Pi)\\
=(\exp(e^{-ad(q)}u')\exp(u),\exp(q')\exp(q),\Pi)\,.
\end{align*}
$$
\xymatrix@C=13pt@R=-2pt@M=6pt{
& \ar@{=>}[dddd]^{e^u} && \ar@{=>}[dddd]^{e^{u'}} &&&&&&& \\
&&&&&&& \ar@{=>}[dd]^{e^{Ad(e^{-q})u'}e^u} &&& \\
\Pi \ar@/^2pc/ [rr]^{e^q} \ar@/_2pc/ [rr]_{e^qe^{d_\Pi u}}
& & e^q\cdot\Pi \ar@/^2pc/ [rr]^{e^{q'}} \ar@/_2pc/ [rr]_{e^{q'}e^{d_\Pi u'}} 
& & e^{q'}\cdot(e^q\cdot\Pi) \ar@{~>}[rr]
& & \Pi\ar@/^2pc/ [rrrr]^{e^{q'}e^q} \ar@/_2pc/ [rrrr]_{e^{q'}e^{d_\Pi u'}e^qe^{d_\Pi u}} 
&&&& (e^{q'}e^q)\cdot\Pi \\ 
&&&&&&&&&& \\
&&&&&&&&&&}
$$
$$
\texttt{Figure 2.~Horizontal composition}
$$
~\\

\subsection{On the $\pi_0$ of the total space of a cosimplicial $2$-groupoid}

\subsubsection{The nerve of a $2$-groupoid}

In this paragraph we review the nerve construction for strict 2-groupoids (see \cite{MoSe}).
Let $\mathcal G$ be a 2-groupoid. The {\it nerve} of $\mathcal{G}$, denoted by $N\mathcal{G}$, is the 
simplicial set defined as follows: $0$-simplices of $N\mathcal{G}$ are objects of $\mathcal{G}$, $1$-simplices 
are $1$-arrows in $\mathcal{G}$, $2$-simplices are diagrams of the following form 
\begin{equation}\label{eq-2-simplex}
\xymatrix@C=6pt@R=14pt@M=6pt{
& \bullet \ar[dr]^{a_{23}} \ar@{=>}[d] &  \\
\bullet \ar[rr]_{a_{13}} \ar[ru]^{a_{12}} & & \bullet\,,}
\end{equation}
3-simplices are commutative tetrahedra of the form 
\begin{equation}\label{eq-3-simplex}
\xymatrix@C=3pt@R=12pt@M=6pt{
&&& \bullet &&& \\
&&&&  && \\
&&& \bullet \ar[uu]_{t_{24}} \ar@{=>}[ul]|<{t_{124}} \ar@{=>}[d]|{t_{123}} \ar[drrr]_{t_{23}} & & & \\
\bullet \ar[rrrrrr]_{t_{13}} \ar[urrr]^(0.75){t_{12}} \ar[rrruuu]^{t_{14}} && 
\ar@{:>}[ul]|{t_{134}} &&& & \ar@{=>}[uulll]|{t_{234}} 
\bullet \ar[llluuu]_{t_{34}}\,,}
\end{equation}
and for $n\geq 3$ an $n$-simplex of $N\mathcal{G}$ is an $n$-simplex such that each of its sub-3-simplices 
is a commutative tetrahedron as above. 

$N\mathcal{G}$ is fibrant and such that $\pi_0(N\mathcal{G})$ is the quotient set of $\mathcal{G}$, 
$\pi_1(N\mathcal{G},x)$ is the group of $1$-automorphisms of $x$, $\pi_2(N\mathcal{G},x)$ is the 
group of $2$-automorphisms of ${\rm id}_x$, and $\pi_n(N\mathcal{G},x)=0$ for $n\geq3$. 

\subsubsection{Path components of the total space of a cosimplicial $2$-groupoid}

Let $\mathcal{G}^\bul$ be a fibrant cosimplicial $2$-groupoid, meaning that $N\mathcal{G}^\bul$ is fibrant 
as a cosimplicial simplicial set. In this paragraph we describe the set $\pi_0({\rm Tot}(N\mathcal{G}^\bul))$. 

It follows from the nerve construction of the previous paragraph that $\big({\rm Tot}(N\mathcal{G}^\bul)\big)_0$ 
is the set of $4$-tuples $(m,g,a,t)$, where 
\begin{itemize}
\item $m$ is an object in $\mathcal G^0$, 
\item $g$ is a $1$-arrow in $\mathcal{G}^1$ with source $\delta_2m$ and target 
$\delta_1m$, 
\item $a$ is a diagram of the form \eqref{eq-2-simplex} in $\mathcal{G}^2$ such 
that $a_{12}=\delta_3g$, $a_{13}=\delta_2g$ and $a_{23}=\delta_1g$, 
\item $t$ is a tetrahedron \eqref{eq-3-simplex} in $\mathcal{G}^3$ such that 
$t_{123}=\delta_4a$, $t_{134}=\delta_2a$, $t_{124}=\delta_3a$ and 
$t_{234}=\delta_1a$. 
\end{itemize}
Then $\pi_0({\rm Tot}(N\mathcal{G}^\bul))=\big({\rm Tot}(N\mathcal{G}^\bul)\big)_0/\sim$, 
where two $0$-simplices $(m,g,a,t)$ and $(m',g',a',t')$ are equivalent through $\sim$ 
if there exists a triple $(\gamma,\alpha,\tau)$ such that 
\begin{itemize}
\item $\gamma$ is a $1$-arrow in $\mathcal G^0$ with source $g$ and target $g'$, 
\item $\al$ is a diagram of the form 
$$
\xymatrix{
\bul \ar[d]_{\delta_2\gamma} \ar[rr]^{g} & \ar@{=>}[d] & \bul \ar[d]^{\delta_1\gamma} \\
\bul \ar[rr]_{g'} & & \bul\,,}
$$
\item $\tau$ is a commutative diagram of the following form 
$$
\xymatrix@C=6pt@R=1pt@M=5pt{
& & & \bullet \ar[dddd] \ar[dddrr] & & \\
& & \ar@{:>}[dd]^{a} & & & \\
& \ar@{=>}[dddd]^{\tau_{12}} & & & \ar@{=>}[ddd]^{\tau_{23}} &  \\
\bullet \ar[dddd] \ar@{.>}[rrrrr] \ar[rrruuu] & & \ar@{:>}[dddd]^(0.65){\tau_{13}} & & & \bullet \ar[dddd] \\
& & & \bullet \ar@{=>}[ddd]^(0.60){a'} \ar[dddrr] & & \\
& & & & & \\
& & & & & \\
\bullet \ar[rrrrr] \ar[rrruuu] & & & & & \bullet\,,}
$$
where $\tau_{12}=\delta_3\gamma$, $\tau_{13}=\delta_2\gamma$ and $\tau_{23}=\delta_1\gamma$. 
\end{itemize}

\subsection{The Maurer-Cartan functor of a quantum type cosimplicial DGLA}

Let $\mg^\bul$ be a fibrant quantum type cosimplicial DGLA. Then for any pronilpotent 
commutative algebra $\mm$, we can now explicitely describe $\underline{MC}_{\mg^\bul}(\mm)$: it is 
the quotient of the set of {\it weak Maurer-Cartan elements} by {\it weak gauge equivalences}, that we define 
in the following two paragraphs. 

\subsubsection{Weak Maurer-Cartan elements}

A weak Maurer-Cartan element is a triple $(\Pi,g,a)$ such that 
\begin{itemize}
\item $\Pi$ is a standard Maurer-Cartan element in $\mg^0$, that is to say $\Pi\in(\mg^0\otimes\mm)^{[1]}$ 
satisfies the Maurer-Cartan equation 
$$
d\Pi+\frac12[\Pi,\Pi]=0\,,
$$
\item $g\in\exp\big((\mg^1\otimes\mm)^{[0]}\big)$ is such that 
$$
g\cdot\big(\delta_2\Pi\big)=\delta_1\Pi\,,
$$
\item $a\in\exp\big((\mg^2\otimes\mm)^{[-1]}_{(\delta_3\delta_2)\Pi}\big)$ is such that 
$$
\delta_1g\,\delta_3g=a^{-1}\cdot\big(\delta_2g\big)
$$
and satisfies 
$$
\delta_4a\,\delta_2a=\delta_4\delta_3g(\delta_1a)\,\delta_3a
$$
\end{itemize}
\begin{ex}[Sheaves]\label{E2.1}
\emph{
If $\mg^\bullet$ is the quantum cosimplicial DGLA naturally associated to a sheaf of quantum 
type DGLAs $\mg(-)$ on a topological space $X$ and an open cover $X=\bigcup_{i\in I}U_i$ then 
a weak Maurer-Cartan element is a triple $(\Pi,g,a)$ as follows: 
$\Pi=(\Pi_i)_{0\leq i\leq m}$ with $\Pi_i\in(\mg(U_i)\otimes\mm)^{[1]}$ satisfying the 
Maurer-Cartan equation for any $i$, 
$g=(g_{ij})_{0\leq i<j\leq m}$ with $g_{ij}\in\exp\big((\mg(U_{ij})\otimes\mm)^{[0]}\big)$ 
such that $g_{ij}\cdot\Pi_i=\Pi_j$ on $U_{ij}$ for any $i<j$, 
$a=(a_{ijk})_{0\leq i<j<k\leq m}$ with $a_{ijk}\in\exp\big((\mg(U_{ijk})\otimes\mm)^{[-1]}_{\Pi_i}\big)$
such that 
\begin{equation}
g_{jk}g_{ij}=a_{ijk}^{-1}\cdot g_{ik}
\end{equation}
on $U_{ijk}$ for any $i<j<k$, 
with the additional condition that 
\begin{equation}
\label{2.1.4}
a_{ijk}a_{ikl}=g_{ij}(a_{jkl})a_{ijl}
\end{equation}
on $U_{ijkl}$ for any $i<j<k<l$. 
}\end{ex}
\begin{ex}[Algebroid stack deformations of sheaves of algebras]\emph{
Let $\cA (-)$ be a sheaf of algebras over a topological space $X$. And let 
$\mg(-)=C(\cA ,\cA )$ be the sheaf of quantum type DGLA given by Hochschild 
cochains. Given a pronilpotent commutative algebra $\mm$ and an open cover $X=\bigcup_iU_i$, 
the corresponding weak Maurer-Cartan elements are (according to the previous example) precisely 
the $\mm$-deformations of $\cA $ as an {\it algebroid stack} (see \cite{K2}). 
}\end{ex}

\subsubsection{Weak gauge equivalences}

A weak (gauge) equivalence between two weak Maurer-Cartan elements $(\Pi,g,a)$ and 
$(\Pi',g',a')$ is a pair $(\gamma,\alpha)$ such that 
\begin{itemize}
\item $\gamma\in\exp(\mg^0\otimes\mm)^{[0]}$ is a gauge equivalence between $\Pi$ and $\Pi'$: 
$$
\Pi'=\gamma\cdot\Pi\,,
$$
\item $\alpha\in\exp(\mg^1\otimes\mm)^{[-1]}_{\delta_2\Pi}$ is such that 
$$
\alpha\cdot\big(\delta_1\gamma\,g\big)=g'\,\delta_2\gamma
$$
and satisfies 
$$
a\,\delta_2\al=\delta_3g(\delta_1\al)\,\delta_3\al\,\delta_3\delta_2\gamma(a')
$$
\end{itemize}
\begin{ex}[Sheaves]\emph{
Let us go back to example \ref{E2.1}. In this context a weak equivalence is a pair $(\gamma,\alpha)$ as follows: 
$\gamma=(\gamma_i)_{0\leq i\leq m}$ with $\gamma_i\in\exp\big(\mg(U_i)\otimes\mm\big)^{[0]}$ being a gauge equivalence between 
$\Pi_i$ and $\Pi'$, and $\alpha=(\alpha_{ij})_{0\leq i<j\leq m}$ with $\alpha_{ij}\in\exp\big((\mg(U_{ij})\otimes\mm)^{[-1]}\big)$ 
such that $\alpha_{ij}\cdot(\gamma_jg_{ij})=g'_{ij}\gamma_i$ on $U_{ij}$ for any $i<j$,  
with the additional condition that 
$$
a_{ijk}\al_{ik}=g_{ij}(\al_{jk})\al_{ij}\gamma_i(a'_{ijk})
$$
on $U_{ijk}$ for any $i<j<k$. }\end{ex}


\subsection{The acyclic case}

In this paragraph we assume  that we are given a cosimplicial quantum type DGLA $\mg^\bullet$ that is 
acyclic as a cosimplicial vector space. Namely $\check{H}^i(\mg^\bullet)=0$ for any $i>0$. 
We have the following: 
\begin{prop}
$\underline{MC}_{\mg^\bullet}(\mm)$ is in one-to-one correspondance with the set of 
usual MC elements up to usual gauge equivalences in the DGLA $\check{H}^0(\mg^\bullet)\otimes\mm$. 
\end{prop}
\begin{proof}
Let us first prove that any weak Maurer-Cartan element $(\Pi,g,a)$ is weakly equivalent to a weak 
Maurer-Cartan element of the form $(\Pi',1,1)$. We do this in using two successive inductions: 
\begin{enumerate}
\item Assume that $a=1~\textrm{mod}~\mm^i$. Then the tetrahedron equation together with 
$\check{H}^2(\mg^\bullet)=0$ implies that $a=1+\check{d}b~\textrm{mod}~\mm^{i+1}$, with 
$b\in(\mg^2\otimes\mm^i)^{[-1]}$. Therefore, applying the weak equivalence $(1,\exp(b))$ 
one obtains a weak Maurer-Cartan element $(\Pi,g',a')$ with $a'=1~\textrm{mod}~\mm^{i+1}$. 
By induction $(\Pi,g,a)$ is weakly equivalent to a Maurer-Cartan element of the form $(\Pi,g',1)$. 
\item Assume that $g'=1~\textrm{mod}~\mm^i$. Then the triangle equation together with 
$\check{H}^1(\mg^\bullet)=0$ implies that $g=1+\check{d}h~\textrm{mod}~\mm^{i+1}$, with 
$h\in(\mg^1\otimes\mm^i)^{[0]}$. Therefore, applying the weak equivalence $(\exp(h),1)$ 
one obtains a weak Maurer-Cartan element $(\Pi',g'',1)$ with $g''=1~\textrm{mod}~\mm^{i+1}$. 
By induction $(\Pi,g',1)$ is weakly equivalent to a Maurer-Cartan element of the form $(\Pi',1,1)$. 
\end{enumerate}

We then observe that if two weak Maurer-Cartan elements $(\Pi,1,1)$ and $(\Pi',1,1)$ are related by a 
weak equivalence $(\gamma,\alpha)$ then the weak equivalence $(\gamma,1)$ also relates them. 

Finally, weak Maurer-Cartan elements of the form $(\Pi,1,1)$ (resp.~weak equivalences of the form 
$(\gamma,1)$) are precisely ususal Maurer-Cartan elements (resp.~usual gauge equivalences) in 
$\check{H}^0(\mg^\bullet)\otimes\mm$. 
\end{proof}

\section{Applications to deformation quantization}

Let $(X,\mO)$ be a topological space equipped with a sheaf of commutative (and associative) unital 
$k$-algebras and assume that $\cL$ is a sheaf of Lie algebroids\footnote{Lie algebroids are also called Lie-Rinehart 
algebras. } over $(X,\mO)$ which is locally free and of constant rank $d\in\N^*$ as an $\mO$-module. 

We also assume that $\mO$ is locally acyclic (so that \v{C}ech resolution is relevant). 

Recall from \cite{C,C2,CDH,CVdB} that there are two sheaves of quantum type DG-Lie algebras $\mg$ and $\mh$ 
associated to $\cL$: $\mg$ is the sheaf of $\cL$-poly-vector fields equipped with zero differential 
and Schouten type bracket $[-,-]$, and $\mh$ is the sheaf of $\cL$-poly-differential operators (a-k-a 
Hochschild cochains associated to $\cL$\footnote{Which are in fact {\it Cartier cochains} \cite{Car} for 
the counital $\mO$-coalgebra $\mathcal U(\cL)$ with values in the trivial bicomodule. }) 
equipped with Hochschild differential $d_H$ and Gerstenhaber bracket $[-,-]_G$ 
(we refer to \cite{C,CDH} for details). 

This is a well-known fact that $\mg$ is the cohomology sheaf of $(\mh,d_H)$, and that the Lie bracket on $\mg$ 
induced from $[-,-]_G$ is precisely $[-,-]$. 

Let us chose an open cover $X=\bigcup_{i\in I}U_i$ of such that $\mO_{|U_i}$ has trivial comology and $\cL_{|U_i}$ 
is a free $\mO_{|U_i}$-module. We then denote by $\mg^\bullet$ and $\mh^\bullet$ the associated cosimplicial 
DGLAs. 

\subsection{(Weak) Poisson structures, deformations and equivalences}

\subsubsection{Definitions}

A {\em Poisson structure} (resp.~{\em formal Poisson structure}) on $\cL$ is a usual Maurer-Cartan element 
in the graded Lie algebra $\check H^0(X,\mg^\bullet)$ (resp.~in $\hbar\check H^0(X,\mg^\bullet)[[\hbar]]$). 

A {\em formal weak Poisson stucture} on $\cL$ is a weak Maurer-Cartan element $(\pi,g,a)$ for $\mg^\bullet$ 
and $\mm=\hbar k[[\hbar]$. One can see that the first order term in the $\hbar$-series $\pi$ is an actual 
Poisson structure. 

A {\em weak formal deformation} (or simply, a {\em weak deformation}) is a weak Maurer-Cartan element for 
$\mh^\bullet$ and $\mm=\hbar k[[\hbar]$, and we denote its class modulo weak gauge equivalences  
by $\underline{\alpha}\in\underline{MC}_{\mh^\bullet}(\hbar k[[\hbar]])$. 

An {\em actual deformation} is a weak Maurer-Cartan element for $\mh^\bullet$ and $\mm=\hbar k[[\hbar]$ 
of the form $(\Pi,g,1)$. An {\em actual equivalence} between two actual deformations is a weak equivalence 
of the form $(\gamma,1)$. 

This terminology is justified by the example of complex manifolds below. 

\subsubsection{Example: complex manifolds}\label{sec-ex-sympl}

Let $X$ be a complex manifold, $\mO=\mO_X$ be the sheaf of holomorphic functions on $X$, and $\cL$ be the sheaf 
of holomorphic vector fields on $X$. Then $\mg$ (resp.~$\mh$) is the DG Lie algebra of usual holomorphic 
poly-vector fields (resp.~poly-differential operators) on $X$. 

Weak deformations then correspond precisely to holomorphic formal deformations of $\mO$ as an algebroid stack
(see the previous section), and NOT as a sheaf of algebras (the latest corresponding to actual deformations). 

Actual deformations correspond to holomorphic formal deformations of $\mO$ as a sheaf of algebras. 

\subsubsection{Quantization problems}

To any actual deformation $(\Pi_\hbar,g_\hbar)$ on $\cL$ one can associate canonically a Poisson structure 
by taking the skew-symmetrization of the first order term in the $\hbar$-series $\Pi_\hbar$. 
The {\em strong quantization problem} is then as follows: 
\begin{pb}\label{pb1}
Let $\pi$ be a Poisson structure on $\cL$. Does there exists an actual deformation with associated Poisson structure 
being $\pi$ ?
\end{pb}
If it exists, such an actual deformation is called an {\em actual quantization} of $\pi$. \\

In full generality the answer to this problem is NO. This leads us to 
formulate a weaker version of this problem. 
As above, to any weak deformation $(\Pi_\hbar,g_\hbar,a_\hbar)$ on $\cL$ one can associate canonically 
a Poisson structure. The {\em weak quantization problem} is then as follows: 
\begin{pb}\label{pb2}
Let $\pi$ be a Poisson structure on $\cL$. Does there exists a weak deformation with associated Poisson structure 
being $\pi$ ?
\end{pb}
If it exists, such a weak deformation is called a {\em weak quantization} of $\pi$. 

\subsection{Existence and classification of weak quantizations}

\subsubsection{A formality theorem}

There exists several results \cite{Ye,CDH,VdB,CVdB} about the extension of Kontsevich formality theorem for algebraic 
varietes and/or complex manifolds. The one we will use in the paper is taken from \cite{CVdB} (Theorem 6.4.1) that we 
translate as follows in the language we are using here: 
\begin{thm}[\cite{CVdB}]\label{formality}
The sheaves of DG-Lie algebras $\mg$ and $\mh$ are weakly equivalent. 
\end{thm}

\begin{rmq}
The graded vector subspace $\tilde\mh\subset\mh$ defined by 
$$
\tilde\mh^{[k]}:=\big(\ker(\epsilon:\mathcal U(\cL)\to\mO)\big)^{\otimes_\mO k+1}
\subset\mathcal U(\cL)^{\otimes_\mO k+1}=\mh^{[k]}
$$
actually is a DG Lie subalgebra of $\mh$, and the inclusion is a (objectwise) quasi-isomorphism : 
as a counital $\mO$-coalgebra $\mathcal U(\cL)$ is the cofree counital cocommutative and coassociative 
$\mO$-coalgebra (this is PBW Theorem for Lie algebroids), for which this result is a standard fact 
(see e.g.~\cite{CE}, Ch. IX). Thus for any open $U\subset X$, $Del_{\tilde\mh(U)}$ and $Del_{\mh(U)}$ 
are weakly equivalent, and then we can work with $\tilde\mh$ instead of $\mh$. This is more convenient 
since cochains of positive degree in $\tilde\mh$ are vanishing when acting on $k\subset\mO$ in (at least) 
one argument. \footnote{In the case when $\mO=\mO_X$ and $\cL=\mathcal T_X$ are the structure and tangent 
sheaf of a smooth algebraic variety $X$, elements of $\tilde\mh$ are called {\it normalized} poly-differential 
operators in \cite{Ye}. } 
\end{rmq}

\subsubsection{Main result}

As a direct consequence of Theorem \ref{formality} and of discussion in the previous section we have the following result, 
which in particular gives a positive answer to Problem \ref{pb2}: 
\begin{thm}\label{thm-main}
1) Any Poisson structure $\pi$ on $\cL$ admits a weak quantization. \\
\indent 2) For any Poisson structure $\pi$ on $\cL$ there is a one-to-one correspondence 
$$
\frac{\{\textrm{w.P.s.}~(\pi_\hbar,g_\hbar,a_\hbar)~\textrm{s.t.}~\pi_\hbar=\hbar\pi+o(\hbar)\}}{\textrm{weak equivalences}}
\longleftrightarrow
\frac{\{\textrm{weak quantizations of }\pi\}}{\textrm{weak equivalences}}\,.
$$
\end{thm}
This result has first been conjectured (and proved ?) by Leitner and Yekutieli (see \cite{LY}) when 
$\mO=\mO_X$ and $\cL=\mathcal T_X$ are respectively the structure and tangent sheaf of a smooth algebraic 
variety (they speak about {\it twisted} things while we write {\it weak}). 

\subsection{Classification in the complex symplectic case}

In this paragraph $X$ is (the underlying topological space of) a complex manifold and $\mO=\mO_X$ 
is the sheaf of holomorphic functions on it. Thanks to the $\bar{\pa}$-Poincar\'e lemma 
the natural inclusion $\mg^\bullet\hookrightarrow(\cA ^{0,*}(\mg^\bullet),\brd)$ is a weak equivalence. 
Moreover, as $\cA ^{0,*}(\mg)$ is the sheaf of sections of a $C^\infty$-bundle it is acyclic; therefore 
$\underline{MC}_{\mg^\bullet}(\hbar k[[\hbar]])$ is discribed by the set of usual Maurer-Cartan elements 
up to usual gauge equivalences in the DGLA of global sections $\hbar\cA ^{0,*}(\mg)(X)[[\hbar]]$ with 
$\brd$ as differential. We write $\mathcal G=(\cA ^{0,*}(\mg)(X),\brd)$ and introduce the following 
bigrading on it: $\mathcal G^{k,l}=\cA ^{0,l}(\mg^{[k]})(X)$. 

Let us suppose that we are given a Poisson structure $\pi$ on a holomorphic Lie algebroid $\cL$ over 
$(X,\mO)$ which is {\em symplectic}. Namely, viewed as an $\mO$-linear map 
$\pi^\sharp:\cL^\vee={\rm Hom}_{\mO}(\cL,\mO)\to\cL$ it is invertible. Consider 
the bundle $\EOm^*=\wedge^*\cL^\vee$ of \emph{$\cL$-differential forms} and recall (see \cite{CDH}) 
that sections of it (\emph{$\cL$-forms} for short) are endowed with the following 
\emph{$\cL$-de Rham differential}: for any $\cL$-$k$-form $\eta$ and $\cL$-vector 
fields $\sigma_0,\dots,\sigma_k$, 
\begin{eqnarray}
\Edif\eta(\sigma_0,\dots,\sigma_k) & := & 
\sum_i(-1)^{i}\rho(\sigma_i)\eta(\sigma_0,\dots,
\hat\sigma_i,\dots,\sigma_k)\label{E-dif} \\
& & +\sum_{i<j}(-1)^{i+j}\eta([\sigma_i,\sigma_j],
\sigma_0,\dots,\hat\sigma_i,\dots,\hat\sigma_j,\dots,\sigma_k)\,. \nonumber
\end{eqnarray}
Let us denote $J:\mg\to\wedge^{*+1}\cL^\vee$ the inverse map of $\pi^\sharp$ extended by taking iterated 
exterior products. Let us denote $\omega$ the image of $\pi$ through the map $J$. By direct computation 
one gets that $J$ sends the differential $[-,\pi]$ onto the $\cL$-de Rham differential $\Edif$. 
Recall also that if $u$ is a $\cL$-polyvector field then one can define contraction $\io_u$ with $u$ and 
$\cL$-Lie derivative $\Elie_u$ by $u$. These operations are related by the following formulas: 
for $u$ and $v$ of homogeneous degree $k$ and $l$ one has 
$$
\Elie_u=\Edif\circ \io_u+(-1)^k \io_u\circ\Edif\,,
$$
$$
\Elie_u\circ\Elie_v-(-1)^{kl}\Elie_v\circ\Elie_u=\Elie_{[u,v]}\,,
$$
and
$$
\Elie_u\circ \io_v -(-1)^{k(l+1)}\io_v\circ \Elie_u= (-1)^k
\io_{[u,v]}\,.
$$
Let us introduce the bicomplex $\EA ^{*,*}$: 
$\EA ^{k,l}=\cA ^{0,l}(\wedge^k\cL^\vee)$ with differentials $\Edif$ and $\brd$. 
It is naturally equipped with a descending filtration: 
$$
F^p(\EA)=\oplus_{k\geq p}\EA^{k,*}\,.
$$
\begin{thm}\label{prop44}
There is a one-to-one correspondence
$$
\frac{\{\textrm{weak quantizations of}~\pi\}}{\textrm{weak equivalences}}
\longleftrightarrow
\frac1\hbar\omega+H_{tot}^2(F^1(\EA))+\hbar H_{tot}^2(\EA)[[\hbar]]. 
$$
\end{thm}
\begin{proof}
We know from Theorem \ref{thm-main} that the set of weak quantizations of $\pi$ up to weak equivalences 
is in bijection with the set of formal weak Poisson structures $(\pi_\hbar,g_\hbar,a_\hbar)$ such that 
$\pi_\hbar=\hbar\pi+o(\hbar)$ up to weak equivalences. And we have just seen that this set is itself in 
bijection with the set of usual Maurer-Cartan elements $\widetilde\pi_\hbar$ such that 
$\widetilde\pi_\hbar^{1,0}=\hbar\pi+o(\hbar)$ up to usual gauge equivalences in $\hbar\mathcal G[[\hbar]]$. 

Let us write $\pi_\hbar=\widetilde\pi_\hbar^{1,0}$, $q_\hbar=\widetilde\pi_\hbar^{0,1}$ and 
$r_\hbar=\widetilde\pi_\hbar^{-1,2}$. 
The Maurer-Cartan equation reads 
$$
\begin{array}{llll}
(a) & [\pi_\hbar,\pi_\hbar]=0 & (b) & \bar{\pa}(\pi_\hbar)+[q_\hbar,\pi_\hbar]=0 \\ 
(c) & \bar{\pa}(q_\hbar)+[\pi_\hbar,r_\hbar]+\frac{1}{2}[q_\hbar,q_\hbar]=0 & 
(d) & \bar{\pa}(r_\hbar)+[q_\hbar,r_\hbar]=0\,.
\end{array}
$$

We now define a bundle isomorphism 
\begin{eqnarray*}
T_\C^\vee M & \tilde\longrightarrow & T^{1,0}M\oplus {T^\vee}^{0,1}M \\
\xi+\bar\xi & \longmapsto & \pi_\hbar^\sharp(\xi)-\iota_{q_\hbar}(\xi)+\bar\xi
\end{eqnarray*}
whose inverse extends to a graded (but NOT bigraded) algebra isomorphism 
$$
\widetilde J_\hbar\,:\,\mathcal G[-1]\,\tilde\longrightarrow\,\EA\,.
$$
Then let $\widetilde\omega_\hbar=\frac1\hbar\omega+O(1)$ be the image of $\widetilde\pi_\hbar$ under 
$\widetilde J_\hbar$ and write $\omega_\hbar=\widetilde\omega_\hbar^{2,0}$, $v_\hbar=\widetilde\omega_\hbar^{1,1}$ and 
$u_\hbar=\widetilde\omega_\hbar^{0,2}$. One has $\omega_\hbar=\frac1\hbar\omega+O(1)$, $v_\hbar=O(1)$ and 
$u_\hbar=O(\hbar)$. 
\begin{lem}
1.~$\omega_\hbar=J_\hbar(\pi_\hbar)$ (where $J_\hbar$ is defined as $J$, with $\pi$ 
replaced by $\pi_\hbar$). \\
2.~$q_\hbar=-\pi_\hbar^\sharp(v_\hbar)$. \\
3.~$r_\hbar=u_\hbar-\frac12(\iota_{q_\hbar}^2)(\omega_\hbar)=u_\hbar-\frac12\iota_{q_\hbar}v_\hbar$. 
\end{lem}
\begin{proof}[Proof of the lemma]
The first part is well-known. 

Before proving the second and the third parts observe that for any $q\in\mathcal G^{0,1}$ one has 
$(\iota_q\otimes{\rm id})(\omega)=\frac12\iota_q(\omega)=({\rm id}\otimes\iota_q)(\omega)$ and 
$\iota_q(\omega)=-J(q)$. 

Therefore $q_\hbar=\pi_\hbar^\sharp(v_\hbar)
-(\iota_{q_\hbar}\otimes\pi_\hbar^\sharp+\pi_\hbar^\sharp\otimes\iota_{q_\hbar})(\omega_\hbar)
=\pi_\hbar^\sharp(v_\hbar)+2q_\hbar$, which proves part 2. 
And finally $r_\hbar=u_\hbar-\iota_{q_\hbar}(v_\hbar)+(\iota_{q_\hbar}\otimes\iota_{q_\hbar})(\omega_\hbar)
=u_\hbar-\frac12\iota_{q_\hbar}^2(\omega_\hbar)$. 
\end{proof}
We now prove
\begin{prop}
$(\Edif+\bar{\pa})(\widetilde\omega_\hbar)=0$ if and only if $\widetilde\pi_\hbar$ is a Maurer-Cartan element. 
\end{prop}
\begin{proof}[Proof of the proposition]
First of all, by applying $J_\hbar$ to $(a)$ one sees that it is equivalent to $\Edif(\omega_\hbar)=0$. 

Then $J_\hbar(\brd(\pi_\hbar)+[\pi_\hbar,q_\hbar])=-\brd(\omega_\hbar)+\Edif(J_\hbar(q_\hbar))
=-\brd(\omega_\hbar)-\Edif(v_\hbar)$. Therefore (b) is equivalent to $\brd(\omega_\hbar)+\Edif(v_\hbar)=0$. 

Now assume that (a) and (b) are satisfied. The following lemma tells us that, under this assumption, 
(c) is equivalent to $\brd(v_\hbar)+\Edif(u_\hbar)=0$. 
\begin{lem}
The l.h.s.~of (c) is equal to the $(0,2)$-part of $\widetilde J_\hbar^{-1}(\Edif+\brd)(\widetilde\omega_\hbar)$. 
\end{lem}
\begin{proof}[Proof of the lemma]
The $(0,2)$-part of $\widetilde J_\hbar^{-1}(\Edif+\brd)(\widetilde\omega_\hbar)$ is 
\begin{eqnarray}
&& \pi_\hbar^\sharp(\brd v_\hbar+\Edif u_\hbar)
-(\pi_\hbar^\sharp\otimes\iota_{q_\hbar}+\iota_{q_\hbar}\otimes \pi_\hbar^\sharp)(\brd\omega_\hbar+\Edif v_\hbar) 
\nonumber \\
&& + (\pi_\hbar^\sharp\otimes\iota_{q_\hbar}\otimes\iota_{q_\hbar}
+\iota_{q_\hbar}\otimes \pi_\hbar^\sharp\otimes\iota_{q_\hbar}
+\iota_{q_\hbar}\otimes\iota_{q_\hbar}\otimes \pi_\hbar^\sharp)(\Edif\omega_\hbar) 
\label{eq-bigbig}\\
& = & \pi_\hbar^\sharp(\brd v_\hbar+\Edif u_\hbar)-\pi_\hbar^\sharp\iota_{q_\hbar}(\brd\omega_\hbar+\Edif v_\hbar)
+\frac12\pi_\hbar^\sharp\iota_{q_\hbar}^2(\Edif\omega_\hbar)\,. \nonumber
\end{eqnarray}
We now compute each term separately. First of all 
$$
\pi_\hbar^\sharp(\brd v_\hbar+\Edif u_\hbar)=[\pi_\hbar,u_\hbar]+\pi_\hbar^\sharp(\brd v_\hbar)\,.
$$
Then 
\begin{eqnarray*}
\pi_\hbar^\sharp\iota_{q_\hbar}(\brd\omega_\hbar+\Edif v_\hbar) 
& = & \pi_\hbar^\sharp(\Edif\iota_{q_\hbar}v_\hbar-\Elie_{q_\hbar}v_\hbar
+\brd\iota_{q_\hbar}\omega_\hbar+\iota_{\brd q_\hbar}\omega_\hbar) \\
& = & [\pi_\hbar,\iota_{q_\hbar}^2\omega_\hbar]-\pi_\hbar^\sharp(\Elie_{q_\hbar}v_\hbar)
+\pi_\hbar^\sharp(\brd v_\hbar)-\brd(q_\hbar)\,.
\end{eqnarray*}
Finally
\begin{eqnarray*}
\pi_\hbar^\sharp\iota_{q_\hbar}^2(\Edif\omega_\hbar) 
& = & \pi_\hbar^\sharp\iota_{q_\hbar}(\Edif\iota_{q_\hbar}-\Elie_{q_\hbar})(\omega_\hbar) 
=\pi_\hbar^\sharp(\Edif\iota_{q_\hbar}^2-2\Elie_{q_\hbar}\iota_{q_\hbar}
-\iota_{[q_\hbar,q_\hbar]})(\omega_\hbar) \\
& = & [\pi_\hbar,\iota_{q_\hbar}^2\omega_\hbar]-2\pi_\hbar^\sharp(\Elie_{q_\hbar}v_\hbar)+[q_\hbar,q_\hbar]\,.
\end{eqnarray*}
Therefore the r.h.s.~of \eqref{eq-bigbig} gives 
$$
\brd(q_\hbar)+[\pi_\hbar,u_\hbar-\frac12\iota_{q_\hbar}^2\omega_\hbar]+\frac12[q_\hbar,q_\hbar]
=\brd(q_\hbar)+[\pi_\hbar,r_\hbar]+\frac12[q_\hbar,q_\hbar]\,,
$$
that is precisely the l.h.s.~of (c). The lemma is proved. 
\end{proof}
We assume finally that (a) (b) and (c) are satisfied. The next lemma implies that, under this assumption, 
(d) is equivalent to $\brd(u_\hbar)=0$. 
\begin{lem}
The l.h.s.~of (d) is equal to the $(-1,3)$-part of $\widetilde J_\hbar^{-1}(\Edif+\brd)(\omega_\hbar)$. 
\end{lem}
\begin{proof}[Proof of the lemma]
The $(-1,3)$-part of $\widetilde J_\hbar^{-1}(\Edif+\brd)(\omega_\hbar)$ is 
\begin{eqnarray*}
&& -\iota_{q_\hbar}^{\otimes3}(\Edif\omega_\hbar)+\iota_{q_\hbar}^{\otimes2}(\brd\omega_\hbar+\Edif v_\hbar)
-\iota_{q_\hbar}(\brd v_\hbar+\Edif u_\hbar)+\brd u_\hbar \\
& = & -\frac16\iota_{q_\hbar}^3(\Edif\omega_\hbar)+\frac12\iota_{q_\hbar}^2(\brd\omega_\hbar+\Edif v_\hbar)
-\iota_{q_\hbar}(\brd v_\hbar+\Edif u_\hbar)+\brd u_\hbar\,.
\end{eqnarray*}
As for the previous lemma one computes each term separately. For the reader's convenience we give the results 
without computations: 
\begin{itemize}
\item $\iota_{q_\hbar}^3(\Edif\omega_\hbar)
=-3(\Elie_{q_\hbar}\iota_{q_\hbar}+\iota_{[q_\hbar,q_\hbar]})\iota_{q_\hbar}\omega_\hbar$, 
\item $\iota_{q_\hbar}^2(\brd\omega_\hbar+\Edif v_\hbar)
=(\brd\iota_{q_\hbar}+2\iota_{\brd q_\hbar}-2\Elie_{q_\hbar}\iota_{q_\hbar}
-\iota_{[q_\hbar,q_\hbar]})\iota_{q_\hbar}\omega_\hbar$,
\item $\iota_{q_\hbar}(\brd v_\hbar+\Edif u_\hbar)=(\brd\iota_{q_\hbar}+\iota_{\brd q_\hbar})\iota_{q_\hbar}\omega_\hbar
-\Elie_{q_\hbar}u_\hbar$. 
\end{itemize}
Therefore the $(-1,3)$-part of $\widetilde J_\hbar^{-1}(\Edif+\brd)(\omega_\hbar)$ is 
$$
\brd(u_\hbar-\frac12\iota_{q_\hbar}^2\omega_\hbar)+\Elie_{q_\hbar}(u_\hbar-\frac12\iota_{q_\hbar}^2\omega_\hbar)
=\brd(r_\hbar)+[q_\hbar,r_\hbar]\,,
$$
i.e.~the l.h.s.~of (d). The lemma is proved. 
\end{proof}
This ends the proof of the proposition. 
\end{proof}
Therefore we obtain a map from the set of Maurer-Cartan elements $\widetilde\pi_\hbar$ in $\hbar\mathcal G[[\hbar]]$ 
such that $\widetilde\pi_\hbar^{1,0}=\hbar\pi+o(\hbar)$ to the set of $2$-cocycles 
$\widetilde\omega_\hbar=\frac1\hbar\omega+O(1)$ 
in $(^\cL\cA ^{*,*},\Edif+\bar{\pa})$ with $(0,2)$-part being zero mod $\hbar$. 
This map is obviously bijective. 
Let us now assume that we have another Maurer-Cartan element $\widetilde\pi_\hbar'$ in $\hbar\mathcal G[[\hbar]]$ 
such that $(\widetilde\pi_\hbar')^{1,0}=\hbar\pi+o(\hbar)$, with image under the above map denoted 
$\widetilde\omega_\hbar'$. 
It remains to prove the following: 
\begin{prop}
$\widetilde\pi_\hbar$ and $\widetilde\pi_\hbar'$ are gauge equivalent if and only if 
$\widetilde\omega_\hbar'=\widetilde\omega_\hbar+(\Edif+\bar{\pa})(\widetilde\theta_\hbar)$ with 
$\widetilde\theta_\hbar$ a $1$-cochain who's $(0,1)$-part is zero mod $\hbar$. 
\end{prop}
\begin{proof}[Proof of the proposition]
The gauge equivalence between $\widetilde\pi_\hbar$ and $\widetilde\pi_\hbar'$ can be reformulated as follows: 
$\widetilde\pi_\hbar'=\widetilde\pi_{t|t=1}$ where $\widetilde\pi_t$ is the solution of the differential equation 
\begin{equation}\label{eq-moser}
\frac{d\widetilde\pi_t}{dt}-(\brd\widetilde\alpha_t+[\widetilde\pi_t,\widetilde\alpha_t])=0
\end{equation}
with initial condition $\widetilde\pi_{t|t=0}=\widetilde\pi_\hbar$ and 
$\widetilde\alpha_t\in\exp(\hbar\mathcal G^0[[\hbar]])$. As above we write $\pi_t=\widetilde\pi_t^{1,0}$, 
$q_t=\widetilde\pi_t^{0,1}$ and $r_t=\widetilde\pi_t^{-1,2}$, and we define $\widetilde J_t$ in the same way as 
$\widetilde J_\hbar$. We also write $\widetilde\omega_t=\widetilde{J}_t(\widetilde\pi_t)$, $\omega_t=\widetilde\omega_t^{2,0}$,
$v_t=\widetilde\omega_t^{1,1}$
and $u_t=\widetilde\omega_t^{0,2}$. 
We finally write $\alpha_t=\widetilde\alpha_t^{0,0}$ and $\beta_t=\widetilde\alpha_t^{-1,1}$, and define 
$\widetilde\theta_t=\widetilde{J}_t(\widetilde\alpha_t)=\theta_t+\gamma_t$ ($\theta_t={\widetilde\theta_t}^{1,0}$ and 
$\gamma_t={\widetilde\theta_t}^{0,1}$). So we have
$\alpha_t=\pi_t^\sharp(\theta_t)$ and $\beta_t=
\gamma_t-\iota_{q_t}(\theta_t)$. 
It suffices to prove that the differential equation \eqref{eq-moser} is equivalent to 
\begin{equation}\label{eq-moser2}
\frac{d\widetilde\omega_t}{dt}-(\Edif+\bar{\pa})(\widetilde\theta_t)=0\,.
\end{equation}
First of all the $(1,0)$-part of \eqref{eq-moser} is equivalent 
to the $(2,0)$-part of \eqref{eq-moser2}, i.e.~one has
$$
\dot{\pi_t}-[\pi_t,\alpha_t]=0\Longleftrightarrow\dot{\omega_t}-\Edif(\theta_t)=0\,,
$$ 
which directly follows from the fact that $J_t(\dot{\pi_t})=-\dot{\omega_t}+\Edif(\theta_t)$. 

Let us now assume that the $\dot{\omega_t}-\Edif(\theta_t)=0$.  
\begin{lem}
Under the previous assumption the $(0,1)$-part of \eqref{eq-moser} is equivalent 
to the $(1,1)$-part of \eqref{eq-moser2}.
\end{lem}
\begin{proof}[Proof of the lemma]
Since we assumed that $\dot{\omega_t}-\Edif(\theta_t)=0$ then it suffices to prove that the $(0,1)$-part of 
${\widetilde{J}}_t^{-1}(\dot{\widetilde{\omega_t}}-(\Edif+\bar{\pa})(\widetilde\theta_t))$ 
is equal to the opposite of the $(0,1)$-part of the l.h.s.~of \eqref{eq-moser}. 
The $(0,1)$-part of ${\widetilde{J}}_t^{-1}(\dot{\widetilde{\omega_t}}-(\Edif+\bar{\pa})(\widetilde\theta_t))$ is
\begin{eqnarray*}
&& \pi_t^\sharp(\dot{v_t}-\bar{\pa}\theta_t-\Edif \gamma_t)
-(\pi_t^\sharp \otimes \iota_{q_t}+\iota_{q_t}\otimes\pi_t^\sharp)(\dot{{\omega_t}}-\Edif\theta_t)\\
& = & \pi_t^\sharp(\dot{v_t}-\bar{\pa}\theta_t-\Edif\gamma_t)-\pi_t^\sharp\iota_{q_t}(\dot{{\omega_t}}-\Edif\theta_t)\,.
\end{eqnarray*}
Let us compute the first term of this sum:
\begin{eqnarray*}
\pi_t^\sharp(\dot{v_t}-\bar{\pa}\theta_t-\Edif \gamma_t)
& = & \overbrace{\pi_t^\sharp v_t}^{\cdot}-\dot{\pi_t}^\sharp(v_t)
+\bar{\pa}(\pi_t^\sharp\theta_t)-(\bar{\pa}\pi_t^\sharp)(\theta_t)+[\pi_t,\gamma_t]\\
& = & -\dot{q_t}-\dot{\pi_t}^\sharp(v_t)+\bar{\pa}(\alpha_t)+([q_t,\pi_t]^\sharp)(\theta_t)+[\pi_t,\gamma_t]\,.
\end{eqnarray*}
Then the second term is 
\vskip-1cm
\begin{eqnarray*}
\pi_t^\sharp \iota_{q_t}(\dot{{\omega_t}}-\Edif\theta_t)
& = & \overbrace{\pi_t^\sharp\iota_{q_t}\omega_t}^{\cdot}
-\pi_t^\sharp\iota_{\dot{q_t}}\omega_t-\dot{\pi_t}^\sharp\iota_{q_t}{{\omega_t}}
-\pi_t^\sharp(\Edif\iota_{q_t}\theta_t-\Elie_{q_t}\theta_t)\\
& = & -\dot{q_t}+\dot{q_t}-\dot{\pi_t}^\sharp v_t+ [\pi_t,\iota_{q_t}\theta_t]+\pi_t^\sharp \Elie_{q_t}\theta_t \\
& = & -\dot{\pi_t}^\sharp v_t+[\pi_t,\iota_{q_t}\theta_t]-[q_t,\alpha_t]+([q_t,\pi_t]^\sharp)\theta_t\,.
\end{eqnarray*}
Therefore the
$(0,1)$-part of ${\widetilde{J}}_t^{-1}(\dot{\widetilde{\omega_t}}-(\Edif+\bar{\pa})(\widetilde\theta_t))$ is
$$
-\dot{q_t}+\bar{\pa}(\alpha_t)+[\pi_t,\gamma_t-\iota_{q_t}\theta_t]+[q_t,\alpha_t]\\
=-\dot{q_t}+\bar{\pa}(\alpha_t)+[\pi_t,\beta_t]+[q_t,\alpha_t]\,,
$$
which is minus the $(0,1)$-part of the l.h.s.~of \eqref{eq-moser}.
\end{proof}
To end the proof, assuming \eqref{eq-moser} is true for degrees $(1,0)$ and $(0,1)$ (and \eqref{eq-moser2} is true for degrees $(2,0)$ and $(1,1)$), we have to prove the following lemma:
\begin{lem}
Under the preceding assumptions the $(-1,2)$-part of \eqref{eq-moser} 
is equivalent to the $(0,2)$-part of \eqref{eq-moser2}.
\end{lem}
\begin{proof}[Proof of the lemma]
Thanks to the assumptions we made, it is sufficient to prove that the $(-1,2)$-part of 
${\widetilde{J}}_t^{-1}(\dot{\widetilde{\omega_t}}-(\Edif+\bar{\pa})(\widetilde\theta_t))$ 
is equal to the $(-1,2)$-part of the l.h.s.~of \eqref{eq-moser}. 
The $(-1,2)$-part of ${\widetilde{J}}_t^{-1}(\dot{\widetilde{\omega_t}}-(\Edif+\bar{\pa})(\widetilde\theta_t))$ is
\begin{eqnarray*}
&&\iota_{q_t}^{\otimes 2} (\dot{{\omega_t}}-\Edif\theta_t)
-\iota_{q_t}(\dot{v_t}-\bar{\pa}\theta_t-\Edif\gamma_t)
+\dot{u_t}-\bar{\pa}\gamma_t \\
& =& \frac{1}{2}\iota_{q_t}^{2}(\dot{{\omega_t}}-\Edif\theta_t)
-\iota_{q_t}(\dot{v_t}-\bar{\pa}\theta_t-\Edif\gamma_t)
+\dot{u_t}-\bar{\pa}\gamma_t\,.
\end{eqnarray*}
\end{proof}
Let us compute the first term of this sum:
\begin{eqnarray*}
\frac{1}{2}\iota_{q_t}^{2}(\dot{{\omega_t}}-\Edif\theta_t)
& = & \frac{1}{2}\overbrace{\iota_{q_t}^{2}\omega_t}^{\cdot}-\iota_{\dot{q_t}}\iota_{q_t}\omega_t
-\frac{1}{2}\iota_{q_t}(\Edif\iota_{q_t}-\Elie_{q_t})\theta_t\\
& =&\frac{1}{2}\overbrace{\iota_{q_t}v_t}^{\cdot}-\iota_{\dot{q_t}}v_t
+(\Elie_{q_t}\iota_{q_t}+\frac{1}{2}\iota_{[q_t,q_t]})\theta_t\,.
\end{eqnarray*}
The second term is
\begin{eqnarray*}
\iota_{q_t}(\dot{v_t}-\bar{\pa}\theta_t-\Edif \gamma_t)
& = & \overbrace{\iota_{q_t}v_t}^{\cdot}-\iota_{\dot{q_t}}v_t-\bar{\pa}(\iota_{q_t}\theta_t)
-\iota_{\bar{\pa} q_t}\theta_t+\Elie_{q_t}\gamma_t.
\end{eqnarray*}
Therefore the $(-1,2)$-part of 
${\widetilde{J}}_t^{-1}(\dot{\widetilde{\omega_t}}-(\Edif+\bar{\pa})(\widetilde\theta_t))$ is equal to:
\begin{eqnarray*}
& & -\frac{1}{2}\overbrace{\iota_{q_t}v_t}^{\cdot}+(\Elie_{q_t}\iota_{q_t}+\frac{1}{2}\iota_{[q_t,q_t]})\theta_t
+\bar{\pa}(\iota_{q_t}\theta_t)+\iota_{\bar{\pa} q_t}\theta_t-\Elie_{q_t}\gamma_t
+\dot{u_t}-\bar{\pa}\gamma_t \\
& = & \dot{r_t}-\Elie_{q_t}\beta_t-\bar{\pa} \beta_t+\iota_{\brd q_t+\frac12[q_t,q_t]}\theta_t
=\dot{r_t}-[q_t,\beta_t]-\bar{\pa}\beta_t-\iota_{[\pi_t,r_t]}\theta_t \\
& = &\dot{r_t}-[q_t,\beta_t]-\bar{\pa}\beta_t-[\pi_t^\sharp(\theta_t),r_t]
=\dot{r_t}-[q_t,\beta_t]-\bar{\pa}\beta_t-[\alpha_t,r_t],
\end{eqnarray*}
which is the the $(-1,2)$-part of the l.h.s.~of \eqref{eq-moser}. 
\end{proof}
The theorem is proved. 
\end{proof}
In the case of a complex symplectic manifold (like in paragraph \ref{sec-ex-sympl}) the differential
$\Edif$ is the usual holomorphic differential $\partial$. Therefore (using the $\bar\partial$-Poincare lemma) 
one obtains the following 
\begin{cor}
There is a one-to-one correspondence
$$
\frac{\{\textrm{w.P.s.}~\Pi_\hbar~\textrm{s.t.}~[\Pi_\hbar]=\hbar[\Pi]+o(\hbar)\}}{\textrm{weak equivalences}}
\longleftrightarrow
\frac{1}{\hbar}\omega+F^1\check{H}^2(X,\C)+\hbar\check{H}^2(X,\C)[[\hbar]]\,.
$$
\end{cor}
This classification result is similar to the one in \cite{Po,PS}. 

\subsection{Existence and classification of actual quantizations}

Recall that a weak deformation is an actual one if and only if the $a_{ijk}$'s (see Example \ref{E2.1}) 
are exactly $1$. Let us denote by $\mO^\cL\subset\mO$ the sheaf of subalgebras of $\cL$-invariants. 

\subsubsection{A sufficient condition for the existence}

\begin{prop}
\label{propex}
Assume that the map $\check H^2(X,\mO^\cL)\to\check H^2(X,\mO)$ (which is given by 
$\mO^\cL\hookrightarrow\mO$) is surjective. Then any Poisson structure admits an 
actual quantization. 
\end{prop}
\begin{proof}
We proceed by induction: we prove that for any $n\geq1$ on can build a weak quantization with 
$a=1+O(\hbar^n)$. It is obvious for $n=1$, and for $n=2$ it follows from the fact 
that the starting Poisson structure is not a weak one. Assume now that we have the result for 
$n\geq2$ and let us write $a=1+\hbar^na_n+O(\hbar^{n+1})$. Taking the coefficient of $\hbar^n$ 
in equation \eqref{2.1.4} we get that $\check d(a_n)=0$. By assumption $a_n=\tilde a_n+\check d(\alpha_n)$, 
where $\tilde a_n$ is a \v{C}ech $2$-cocycle with values in $\mO^\cL$. Using the gauge transformation 
$\exp(\hbar^n \alpha_n)$ we get that the quantization is weakly equivalent to a one with 
$a=1+\hbar^n\tilde a_n+O(\hbar^{n+1})$. 
It is an immediate check that replacing $\tilde a_n$ with $0$, equations of Example \ref{E2.1}
are still satisfied. 
So we get the result for $n+1$. 
\end{proof}
\begin{rmq}
Observe that the operation of replacing an invariant element with $0$ is not a weak equivalence. 
Therefore we have not proved that, under the hypothesis of the proposition, any weak quantization is weakly 
equivalent to an actual one. 
\end{rmq}

Let us come back to the example of \S\ref{sec-ex-sympl} and assume that we are given a holomorphic symplectic 
$2$-form. In this case, as $\mO^\cL=\C$, our condition for the existence of an actual (i.e., in this case, holomorphic) 
quantization is the same as in \cite{NT}: surjectivity of $\check H^2(X,\C)\to\check H^2(X,\mO_X)$. Indeed, Nest and 
Tsygan prove that this condition is sufficient for the existence of a Fedosov connection 
$\nabla=\bar{\partial}+\nabla_0+{\rm ad}A+{\rm ad}B$, where $A$ and $B$ are respectively $(0,1)$- and $(1,0)$-forms 
with values in the Weyl algebra bundle of $X$ (cf.~\cite{NT}, Theorem 5.9: in that theorem only surjectivity of 
$\check H^2(X,\C)\to\check H^2(X,\mO_X)$ is used to prove the existence of the connection). Then they prove 
(Theorem 5.6 and following remarks) that taking flat sections of a Fedosov connection one gets the desired 
quantization. 

\subsubsection{A partial classification result}

\begin{prop}
Assume that $\check H^2(X,\mO)=0$ and that the map $\check H^1(X,\mO^\cL)\to\check H^1(X,\mO)$ (which is given by 
$\mO^\cL\hookrightarrow\mO$) is surjective. Then any Poisson structure admits an actual quantization and we have a 
one-to-one correspondence between equivalence classes of weak quantizations and equivalence classes of actual 
quantizations
\end{prop}
\begin{rmq}
We saw that the condition $\check H^2(X,\mO)=0$ implies that any weak deformation is weakly equivalent 
to an actual one. In the same way one can prove that any formal weak Poisson structure is weakly equivalent 
to an actual one. 
\end{rmq}

\begin{proof}
Thanks to the previous remark, we only need to prove that
a weak equivalence can be replaced by a strong one.
Following the same proof as in Proposition \ref{propex} and
the fact that $\check{H}^1(M,\C)\to\check{H}^1(M,\mO_M)$ is surjective, one can 
replace the functions
involved in the definition of the weak equivalence isomorphisms with
constants ones and so $0$ functions.
\end{proof}

Let us again study the example of paragraph \ref{sec-ex-sympl} ($\mO^\cL=\C$) and assume that we 
are given a holomorphic symplectic $2$-form $\omega$.
In this case, existence and classification of actual (i.e. holomorphic) quantizations is the same 
as in \cite{NT}: thanks to the corollary of Proposition \ref{prop44} and the assumption 
$\check{H}^2(X,\mO_X)=0$, equivalence classes of holomorphic quantizations are in one-to-one correspondence 
with $\frac1\hbar\omega+\check{H}^2(X,\C)[[\hbar]]$ (or $\frac1\hbar\omega+H^2(F^1\cA^{*,*}(X))[[\hbar]]$). 
Note that in Nest and Tsygan's construction one really needs the condition $\check{H}^2(X,\mO_X)=0$ otherwise 
the set of preimages of the surjective map $\check{H}^2(X,\C) \to\check{H}^2(X,\mO_X)$ would be an affine space 
(not a vectorial space). So choices of preimages of this map in the construction of a Fedosov cannot be done 
cannonically by chosing element in the space $H^2(F^1\cA^{*,*}(X))$.

\mbox{}

\noindent\footnotesize{\textsc{Departement Mathematik \\
ETH Z\"urich \\
R\"amistrasse 101, 8092 Z\"urich, Switzerland} \\
\emph{E-mail address}: {\bf damien.calaque@math.ethz.ch}, {\bf calaque@math.univ-lyon1.fr} \\

\noindent\textsc{D\'epartement de Math\'ematiques, Universit\'e de Montpellier 2 \\ 
CC 5149, Place Eug\`ene Bataillon \\
F-34095 Montpellier Cedex 5, France} \\
\emph{E-mail address}: {\bf ghalbout@math.univ-montp2.fr}}


\begin{thebibliography}{99}


\bibitem{BGNT}
P.~Bressler, A.~Gorokhovsky, R.~Nest, and B.~Tsygan, 
Deformations of gerbes on smooth manifolds, 
in {\it K-Theory and Noncommutative Geometry}, 349-392, EMS Series of Congress Reports {\bf 2} (Corti$\tilde{\rm n}$as et al Eds.), 2008, EMS Publishing House. 

\bibitem{C}
D.~Calaque, 
{\it Th\'eor\`emes de formalit\'e pour les alg\'ebro\"ides de Lie et quantification des $r$-matrices dynamiques}, 
PhD thesis, Universit\'e Louis Pasteur Strasbourg 1, 2005. 

\bibitem{C2}
D.~Calaque, 
Formality for Lie algebroids, 
Comm. Math. Phys. {\bf257} (2005), no. 3, 563-578. 

\bibitem{CDH}
D.~Calaque, V.~Dolgushev, and G.~Halbout, 
Formality theorem for Hochschild chains in the Lie algebroid setting, 
Journal fur die reine und angewandte Mathematik {\bf612} (2007), 81-127. 

\bibitem{CVdB}
D.~Calaque and M.~Van den Bergh, 
Hochschild cohomology and Atiyah classes, 
Advances in Mathematics {\bf224} (2010), no. 5, 1839-1889. 

%
%
\bibitem{CE}
H. Cartan and S. Eilenberg, 
Homological algebra., 
Princeton Mathematical Series, 19. Princeton, New Jersey : Princeton University Press XV, 390 p., 1956.

\bibitem{Car}
P. Cartier, 
Cohomologie des coalgebras, S\'eminaire Sophus Lie, Expos\'e {\bf5} (1956). 

\bibitem{Ge1}
E.~Getzler, 
A Darboux theorem for Hamiltonian operators in the formal calculus of variations, 
Duke Math. J. {\bf111} (2002), 535-560. 

\bibitem{Ge2}
E.~Getzler, 
Lie theory for nilpotent $L_\infty$-algebras, 
Ann. of Math. {\bf170} (2009), no. 1, 271-301. 

\bibitem{GK}
P.~Goerss and K.~Schemmerhorn, 
Model Categories and Simplicial Methods, 
preprint \texttt{arXiv:math/0609537}. 

%
\bibitem{Hin}
V.~Hinich, 
DG-coalgebras as formal stacks, 
J. Pure Appl. Algebra {\bf162} (2001), no. 2-3, 209-250. 

\bibitem{Hin2}
V.~Hinich, 
Deformation of sheaves of algebras, 
Adv. Math. {\bf195} (2005), no. 1, 102-164. 

%
%
\bibitem{K1}
M.~Kontsevich, 
Deformation quantization of Poisson manifolds, 
Lett. Math. Phys. {\bf 66} (2003), no. 3, 157-216. 

\bibitem{K2}
M.~Kontsevich,
Deformation quantization of algebraic varieties, 
EuroConfeence Moshe Flato 2000, Part III (Dijon), 
Lett. Math. Phys. {\bf 56} (2001), no. 3, 271-294.

%
\bibitem{MoSe}
I.~Moerdijk and J.~Svensson, 
Algebraic classification of equivariant homotopy 2-types, 
J. Pure Appl. Algebra {\bf89} (1993), 187-216. 

\bibitem{NT}
R.~Nest and B.~Tsygan, 
Formal deformations of symplectic Lie algebroids, deformations of holomorphic structures and 
index theorems, Asian J. of Math. {\bf 5} (2001), no. 4, 599-633. 

\bibitem{Po}
P.~Polesello, 
Classification of deformation-quantization algebroids on complex symplectic manifolds, 
Publ. Res. Inst. Math. Sci. {\bf 44} (2008), 725-748. 

\bibitem{PS}
P.~Polesello and P.~Schapira, 
Stacks of quantization-deformation modules on complex symplectic manifolds, 
Int. Math. Res. Notices {\bf 49} (2004), 2637-2664. 

\bibitem{Q1}
D.~Quillen, 
{\it Homotopical algebra}, 
Lectures Notes in Math. {\bf43}, Springer-Verlag, Berlin-Heidelberg-New York, 1967. 

\bibitem{Q2}
D.~Quillen, 
Rational homotopy theory, 
Ann. of Math. (2) {\bf90} (1969), 205-295. 

\bibitem{Ree}
C.L.~Reedy, 
Homotopy theory of model categories, 
preprint (1973) (\url{http://www-math.mit.edu/~psh/#Reedy}). 

\bibitem{VdB}
M.~Van den Bergh, 
On global deformation quantization in the algebraic case, 
Journal of Algebra {\bf315} (2007), no. 1, 326-395. 

\bibitem{Ye}
A.~Yekutieli, 
Deformation Quantization in Algebraic Geometry, 
Adv. Math. {\bf198} (2005), no. 1, 383-432. 

\bibitem{LY}
A.~Yekutieli, 
Twisted Deformation Quantization of Algebraic Varieties, 
preprint \url{arXiv:0905.0488}. 


\end{thebibliography}
\end{document}